\documentclass[a4paper,11pt]{article}
\usepackage{color,xcolor,ucs}
\usepackage[top=0.2in, bottom=0.43in, left = 0.65in, right = 0.65in]{geometry}
\usepackage[linkcolor=black,colorlinks=true,urlcolor=blue]{hyperref}
\usepackage{mathtools}

\usepackage{color,xcolor,ucs}
\usepackage{mathtools}   \usepackage{tikz} 
\usepackage{ amssymb }
\usepackage{extarrows} 
\usepackage{pgf,tikz}
\usepackage{float}
\usetikzlibrary{positioning}
\usetikzlibrary{shapes.geometric}
\usetikzlibrary{shapes.misc}
\usetikzlibrary{arrows}
\usepackage{caption}
\usepackage{mathrsfs}
\usetikzlibrary{arrows,shapes,automata,backgrounds,petri,positioning}
\usetikzlibrary{decorations.pathmorphing}
\usetikzlibrary{decorations.shapes}
\usetikzlibrary{decorations.text}
\usetikzlibrary{decorations.fractals}
\usetikzlibrary{decorations.footprints}
\usetikzlibrary{shadows}
\usetikzlibrary{calc}
\usetikzlibrary{spy}
\usepackage{amsmath}
\usepackage{array}
\usepackage{ amssymb }
\usepackage{braket}
\usepackage{qcircuit}
\usepackage{soul}
\usepackage{braket} 
\usepackage{relsize}

\usepackage{amsmath}
\usepackage{ amssymb }
\usepackage{braket}
\usepackage{qcircuit}
\usepackage{soul}
\usepackage{braket}

\title{{\LARGE Correlation length lower bound for the random-field Potts model with the greedy lattice animal}}
\author{Pete Rigas}
\date{}

\begin{document}

\maketitle

\begin{abstract}
   Motivated by recent developments over the past few years in the study of the correlation length of the random-field Ising model due to Ding and Wirth in a paper first available in 2020, we pursue one natural direction of research that the authors propose is of interest, namely in confirming that the same scaling for the correlation length for the random-field Ising model equals that of the random-field Potts model. To demonstrate that the $\frac{4}{3}$ emergence in the correlation length scaling for the random-field Potts model coincides with that correlation length scaling for the random-field Ising model, we refer to arguments due to Talagrand, in which an upper bound on the greedy lattice animal readily provides a corresponding lower bound on the correlation length. Contributions from the $\frac{4}{3}$ exponent appearing in the correlation length scaling for the random-field Ising, and random-field Potts models alike are reminiscent of $\frac{4}{3}$ exponents in upper bounds taken under the Liouville quantum gravity metric in high temperature. 
\end{abstract}

\section{Introduction}

\subsection{Overview}

\noindent The random-field Ising model, as a variant of the classic Ising model, has received attention for behavior in three-dimensions at zero temperature {\color{blue}[2]}, lower critical dimensions {\color{blue}[10,11]}, universality in three-dimension permitting for exact computation of all critical exponents {\color{blue}[8]}, correlation inequalities {\color{blue}[4]}, algorithmic complexity and connections to phase transitions {\color{blue}[1]}, long range order {\color{blue}[6]}, replica symmetry breaking {\color{blue}[16]}, supersymmetry {\color{blue}[23]}, thermal ground states {\color{blue}[25]}, and other related developments {\color{blue}[7,10,15,17,20]}. As a generalization of the Ising model, the Potts model instead has spins that take on an arbitrary, strictly positive, number of $q$ values, with its randomized-field version also being subject of interest, particularly relating to characterizations of dipolar-like interactions {\color{blue}[3]}, a phase diagram in three dimensions {\color{blue}[9]}, with related characterizations of phase transitions {\color{blue}[24]}, critical behavior for the three-state, random-field, model {\color{blue}[12]}, approximation of ground states from graph cuts {\color{blue}[13]}, formulation of algorithms for quasi-exact approximation of ground states {\color{blue}[14]}, modeling of polar domain interactions of Magnesium compounds {\color{blue}[18]}, and simulating behaviors via Monte Carlo {\color{blue}[19]}. In {\color{blue}[5]}, the authors make use of arguments provided in {\color{blue}[21,22]} to obtain upper and lower bounds for the correlation length of the random-field Ising model. To extend their arguments to other models of interest in Statistical Mechanics, specifically the random-field Potts model which they remark is expected to have the same correlation length scaling, we introduce the quantities used throughout the argument for obtaining the lower bound on the correlation length.

\subsection{Paper organization}

\noindent We provide an overview of the random-field Potts model (RFPM), and then describe the extension of arguments presented in {\color{blue}[5]} for obtaining the same correlation length scaling obtained for the random-field Ising model.

\section{Obtaining the correlation length lower bound}

\noindent In the second section, we introduce the model, and obtain the desired correlation length lower bound in the next section.

\subsection{RFPM objects}

\noindent In the following, for $N \geq 1$, denote $\Lambda_N \equiv \big\{ v \in \textbf{Z}^2 : |v|_{\infty} \leq N \big\}$ as the box of side length $2N$ centered about the origin, for $\textbf{Z}^2 \equiv G = \big( V , E \big)$.

\bigskip

\noindent \textbf{Definition} \textit{1} (\textit{q-state, Random-field Potts model Hamiltonian}). Introduce, for nearest-neighbors $u , v \in V \big( \textbf{Z}^2 \big)$, 

\begin{align*}
 \mathcal{H}^{\mathrm{RFPM}, \mathrm{w} \backslash \mathrm{f} } \big( s , \Lambda_N ,      q  , h   \big) \equiv   \mathcal{H}^{\mathrm{w}  \backslash \mathrm{f}} \big(   s , \Lambda_N ,      q  ,h  \big)   \equiv      - \beta \bigg( \underset{i ,j \in \Lambda_N}{\underset{i \sim j}
{\sum}}    \delta_{s_i , s_j}       + 
 \underset{i \in \Lambda_N , j \notin \Lambda_N}{\underset{i \sim j }{\sum}}  
 \overset{1-q}{\underset{\alpha=0}{\sum}}   \epsilon h^{\alpha}_i \delta_{s_i,\alpha}       \bigg)  \text{ } \text{ , } 
\end{align*}

\noindent corresponding to the Hamiltonian of the model under the spin space $\big\{ 0 , \cdots , q-1 \big\}^{\Lambda_N}$ with wired, or free, boundary conditions, for inverse temperature $\beta > 0$, $q, \epsilon >0$, and Kronecker delta, with,

\begin{align*}
      \delta_{s_i , s_j} \equiv 1 \Longleftrightarrow s_i \equiv s_j   \text{ } \text{ , }  \\ \delta_{s_i,\alpha} \equiv 1 \Longleftrightarrow s_i \equiv \alpha \in \big\{ 0 , \cdots , q-1 \big\}          \text{ } \text{ , } 
\end{align*}

\noindent for the quenched field,

\begin{align*}
h_i \equiv \underset{\alpha \in \{ 0 , \cdots , q-1 \}}{\bigcup}       h^{\alpha}_i   \text{ } \text{ , }  
\end{align*}

\noindent given by the collection,

\begin{align*}
  h \equiv \underset{i \in \Lambda_N}{\bigcup} h_i \equiv  \underset{\alpha \in \{ 0 , \cdots , q-1 \}}{\underset{i \in \Lambda_N}{\bigcup}}  
h^{\alpha}_i    \equiv \big\{ \forall  i \in \Lambda_N  , \exists  \text{ } {0 \leq \alpha \leq 1-q} \text{ }  : \text{ }      h^{\alpha}_i \sim \textbf{N} \big(     0 , \epsilon^{-2} \big)      \big\} \text{ } \text{ , }
\end{align*}

\bigskip

\noindent \textbf{Definition} \textit{2} (\textit{Random-field Potts model probability measure, and expectation, from the Hamiltonian}). Introduce,

\begin{align*}
    \textbf{P}^{\mathrm{w} \backslash \mathrm{f}}_{\beta, \Lambda_N, h} \big[            s      \big]   \equiv \textbf{P}^{\mathrm{w} \backslash \mathrm{f}}_{\Lambda_N} \big[            s      \big] \equiv   \mathcal{Z}^{-1} \mathrm{exp} \big(  - \beta \mathcal{H}^{\mathrm{w}  \backslash \mathrm{f}} \big( s  , \Lambda_N , q      , h    \big)    \big)                      \text{ } \text{ , }
\end{align*}

\noindent corresponding to the probability measure of sampling the spin configuration $s$ from the random-field Potts model, at inverse temperature $\beta >0$, with partition function $\mathcal{Z}$ so that $\textbf{P}_{\Lambda_N} \big[ \cdot \big]$ is a probability measure. From this base probability measure, we denote the expectation of RFPM with, 

\begin{align*}
    < \cdot >^{\mathrm{w}}   \equiv   < \cdot >_{\textbf{P}^{\mathrm{w}}_{\beta , \Lambda_N , h }}    \text{ } \text{ , } 
\end{align*}

\noindent corresponding to the wired boundary condition expectation, and,

\begin{align*}
    < \cdot >^{\mathrm{f}}   \equiv   < \cdot >_{\textbf{P}^{\mathrm{f}}_{\beta , \Lambda_N , h }}   \text{ } \text{ , } 
\end{align*}

\noindent corresponding to the free boundary condition expectation.

\noindent \textbf{Definition} \textit{3} (\textit{probability measure and expectation for the quenched field}). From the quantity $h$ introduced above with \textbf{Definition} \textit{1}, denote the probability measure, and expectation, with respect to $h$ with $\mathscr{P}$ and $\mathscr{E}$, respectively.

\bigskip

\noindent \textbf{Definition} \textit{4} (\textit{spontaneous magnetization of the random-field Potts model}). Introduce,

\begin{align*}
 \mathscr{M}_{\beta , \Lambda_N , h} \equiv   \mathscr{M} \equiv \frac{1}{2}   \mathscr{E} \big[  <     s^{\mathrm{w}}_{0}    >^{\mathrm{w}} - <   s^{\mathrm{f}}_0    >^{\mathrm{f}}        \big]  \text{ } \text{ , }
\end{align*}

\noindent corresponding to the RFPM spontaneous magnetization, for spins $s_0$ located at the origin.

\bigskip

\noindent \textbf{Definition} \textit{5} (\textit{correlation length}). Introduce, for some $\mathscr{L} \in \textbf{R}$,

\begin{align*}
\mathscr{L} \big(           \beta , N  , \mathcal{M} , h   \big) \equiv   \mathscr{L} \equiv     \underset{N > 0}{\mathrm{min}} \big\{ \forall   N  , \exists \mathcal{M} \in \big( 0 ,1 \big) :   \mathscr{M}_{\beta , \Lambda_N , h}  \leq \mathcal{M}    \big\}      \text{ } \text{ , } 
\end{align*}

\noindent corresponding to the RFPM \textit{correlation length}.

\bigskip

\noindent \textbf{Theorem} \textit{1} (\textit{direct application of greedy lattice animal arguments from Theorem 1.1 of {\color{blue}[5]} \textit{for obtaining lower bounds on the RPMF correlation length}}). For every $\mathcal{M} \in ( 0 ,1)$, there exists suitable $\mathcal{C} \equiv \mathcal{C} \big( \mathcal{M} \big) > 0$ so that $\mathscr{L} \big( \infty, N , \mathcal{M} , h \big) \equiv \mathscr{L} \geq \mathrm{exp} \big( \mathcal{C}^{-1} \epsilon^{\frac{4}{3}}  \big)$, at $\beta = \infty$.

\bigskip

\noindent To prove \textbf{Theorem} \textit{1}, below we specify additional objects.

\bigskip

\noindent \textbf{Definition} \textit{5} (\textit{greedy lattice animal}). Introduce,

\begin{align*}
  \mathcal{G}_N       \equiv \underset{A \in \mathcal{A}_N}{ \mathrm{max}}   \frac{\underset{0 \leq \alpha \leq 1-q}{\underset{i \in A}{\mathlarger{\sum}}}   h^{\alpha}_i   }{| \partial A |}         \text{ } \text{ , } 
\end{align*}

\noindent corresponding to the \textit{greedy lattice animal}, for the collection of all simply connected subsets of $\Lambda_N$,

\begin{align*}
    \mathcal{A}_N \equiv \big\{ 
 \Lambda_N  \supsetneq  A \neq \emptyset : A \cap \{ 0 \} \neq \emptyset \big\}     \text{ } \text{ , } 
\end{align*}

\noindent in addition to the related object,

\begin{align*}
  \mathcal{G}_N \supsetneq  \mathscr{G}_N  \equiv  \underset{A \in \mathscr{A}_N}{\mathrm{max}}   \frac{\underset{0 \leq \alpha \leq 1-q}{\underset{i \in A}{\mathlarger{\sum}}}   h^{\alpha}_i   }{| \partial A |}       \text{ } \text{ , } 
\end{align*}

\noindent corresponding to the \textit{greedy lattice animal} for $ \mathscr{A}_N \subsetneq \mathcal{A}_N$

\bigskip

\noindent From the \textit{greedy lattice animals} $\mathcal{G}_N$ and $\mathscr{G}_N$ provided above, the following result upper bounds the expected value of the \textit{greedy lattice animal} which provides the desired lower bound for the correlation length provided in \textbf{Theorem} \textit{1}.

\bigskip

\noindent \textbf{Theorem} \textit{2} (\textit{upper bounds on the two greedy lattice animals, adaptation of Theorem 1.6 from {\color{blue}[5]}}). There exists strictly positive $C$ such that,

\begin{align*}
         \mathscr{E}       \big[ \mathscr{G}_N \big]  \leq      \mathscr{E} \big[ \mathcal{G}_N \big]  \leq C \big( \mathrm{log} N \big)^{\frac{3}{4}}  \text{ } \text{ , } 
\end{align*}

\noindent for all $N \geq 3$.

\bigskip

\noindent In the next section we turn towards the proofs of \textbf{Theorem} \textit{1}, and of \textbf{Theorem} \textit{2}.

\section{From upper bounds on the greedy lattice animal to lower bounds on the correlation length}

\noindent In the following arguments, we make use of the following procedure for constructing polygons.

\bigskip

\noindent \textbf{Algorithm} (\textit{constructing polygons for obtaining the $\frac{4}{3}$ exponent}). Fix $N \geq \mathrm{exp} \big( C \epsilon^{-\frac{4}{3}} \big)$, and denote,

\begin{align*}
    \underset{0 \leq \alpha \leq 1-q}{\underset{i \in A}{\mathlarger{\sum}}}   h^{\alpha}_i  \equiv   w \big( A \big)      \text{ } \text{ . } 
\end{align*}

\noindent To construct a polygon $P \subset [-N,N]^2$, 

\begin{itemize}
\item[$\bullet$] \underline{Step 1} (\textit{initialization}). Set $P_1 \equiv [ -\frac{N}{2} , \frac{N}{2} ]^2$. 

\item[$\bullet$] \underline{Step 2} (\textit{constructing $P_2$ from $P_1$}). For each side $S$ of $P_1$, construct an isosceles triangle $T_S$, with length $\frac{l(S)}{2}$, and height $\frac{\epsilon^{\frac{2}{3}} l(S)}{8}$. With probability $\frac{1}{2}$, if $w \big( T_S  \big) > 0$, append $T_s$ to $P_1$. Otherwise, if $w \big( T_S  \big) < 0$, partition $S$ into four sides of equal length, from which 
the next polygon $P_2$ is obtained from $P_1$ by applying the procedure described earlier in this step for obtaining triangles of prescribed lengths and heights. 

\item[$\bullet$] \underline{Step 3} (\textit{iterating the procedure to obtain the polygon $P_{k+1}$ from the polygon $P_{k}$}). Given $P_k$, construct $P_{k+1}$ by applying \underline{Step 2} repeatedly.

\end{itemize}

\bigskip

\noindent We also make use of the following proposition, appearing in \textit{Theorem 4.4.2}, of {\color{blue}[21]}, with related arguments in {\color{blue}[22]}.

\bigskip

\noindent \textbf{Proposition} \textit{1} (\textit{adaptation of the Talagrand result for the random-field Potts model greedy lattice animal}). There exists strictly positive $C_1$, such that for $N > N^{\prime} >1$,

\begin{align*}
   \mathscr{P} \bigg[      \underset{B \in \mathcal{B}_N}{\mathrm{max}}        \frac{\underset{0 \leq \alpha \leq 1-q}{\underset{i \in A}{\mathlarger{\sum}}}   h^{\alpha}_i   }{| \partial B |}   > C_1 \big( \mathrm{log} N \big)^{\frac{3}{4}} + u       \bigg]   \leq \mathrm{exp} \big( - \frac{u^2}{2} \big)              \text{ } \text{ , } 
\end{align*}

\noindent for the collection of simply connected \textit{greedy lattice animals},

\begin{align*}
 \mathcal{B}_N    \equiv  \underset{N^{\prime} > 0 : \Lambda_{N^{\prime}} \subset \Lambda}{\bigcup}   \mathcal{B}_{N^{\prime}}                 \text{ } \text{ , } 
\end{align*}

\noindent $\forall u \geq 0$.

\bigskip

\noindent \textit{Proof of Proposition 1}. Follow the arguments of \textit{Theorem 4.4.2} of {\color{blue}[21]}, or, for more explicit arguments, \textit{4} of {\color{blue}[5]}.

\bigskip

\noindent \textit{Proof of Theorem 1, and of Theorem 2}. With the collection of \textit{greedy lattice animals} $\mathcal{B}_N$, it suffices to demonstrate,

\begin{align*}
   \underset{B \in \mathcal{B}_N}{\mathrm{max}}        \frac{\underset{0 \leq \alpha \leq 1-q}{\underset{i \in A}{\mathlarger{\sum}}}   h^{\alpha}_i   }{| \partial B |}   
 \equiv    \underset{B^{\prime} \in \mathcal{B}_N}{\mathrm{max}}        \frac{\underset{0 \leq \alpha \leq 1-q}{\underset{i \in A}{\mathlarger{\sum}}}   h^{\alpha}_i   }{| \partial B^{\prime} |}    \text{ } \text{ , } 
\end{align*}

\noindent given $B , B^{\prime} \in \mathcal{B}_N$. From $B$, denote, for the set of vertices $V \big( B \big) \subsetneq V$,

\begin{align*}
\widetilde{B} \equiv \big\{v_1 \in  V \big( B \big)   :  v_1 \cap \mathcal{I} \big( B \big) \neq \emptyset   \big\} \text{ } \text{ , } 
\end{align*}

\noindent corresponding to the collection of vertices that have nonempty intersection with the vertices enclosed within $B$, for,

\begin{align*}
\mathcal{I} \big( B \big)   \equiv  \mathcal{B}^c \backslash \mathcal{B} \text{ } \text{ . } 
\end{align*}

\noindent Under this assumption on vertices from $\widetilde{B}$,

\begin{align*}
\mathscr{P} \big[       B  \longleftrightarrow    \infty   \big] \equiv 0 \text{ } \text{ . } 
\end{align*}

\noindent Before concluding the argument, additionally denote, $\mathfrak{CC} \big( \widetilde{B} \backslash \widetilde{B} \big)$, corresponding to the connected components of $\widetilde{B} \backslash \widetilde{B}$, from which $\partial B$ can be expressed with,

\begin{align*}
   \partial B =        \partial \widetilde{B}  \text{ } \dot{\cup} \text{ }   \big( \underset{1 \le i \leq k}{\bigcup} \partial B_i  \big)   \text{ } \text{ , } 
\end{align*}

\noindent as a result of the fact that $B_1 , \cdots , B_k$ are simply connected. Hence,

\begin{align*}
     \frac{\underset{0 \leq \alpha \leq 1-q}{\underset{i \in A}{\mathlarger{\sum}}}   h^{\alpha}_i   }{| \partial B |}  \leq     \frac{1}{| \partial B |} \underset{1 \leq i \leq k}{\sum}      \big|   \underset{v \in B_i}{\sum}  h^{\alpha}_i     \big|   \equiv  \frac{1}{| \partial B |} \underset{1 \leq i \leq k}{\sum}      \frac{| \partial B_i |}{| \partial B_i |}       \big|   \underset{v \in B_i}{\sum}  h^{\alpha}_i     \big|  \leq  \frac{1}{| \partial B |}  {\mathrm{max}}        \frac{\underset{0 \leq \alpha \leq 1-q}{\underset{i \in A}{\mathlarger{\sum}}}   h^{\alpha}_i   }{| \partial B |} \leq {\mathrm{max}}        \frac{\underset{0 \leq \alpha \leq 1-q}{\underset{i \in A}{\mathlarger{\sum}}}   h^{\alpha}_i   }{| \partial B |}   \text{ } \text{ , } 
\end{align*}

\noindent because,

\begin{align*}
   \frac{1}{| \partial B|} < 1     \text{ } \text{ , } 
\end{align*}

\noindent and also because,

\begin{align*}
 \underset{i \in B}{\sum} \frac{| \partial B_i|}{| \partial B|}  = 1    \text{ } \text{ , } 
\end{align*}

\noindent from which it follows that,

\begin{align*}
    \underset{B \in \mathcal{B}_N}{\mathrm{max}}        \frac{\underset{0 \leq \alpha \leq 1-q}{\underset{i \in A}{\mathlarger{\sum}}}   h^{\alpha}_i   }{| \partial B |}      \overset{(\textbf{Proposition} \text{ } \textit{1})}{\sim} \big( \mathrm{log}N \big)^{\frac{3}{4}}\text{ } \text{ , } 
\end{align*}

\noindent hence providing the desired upper bound for the \textit{greedy lattice animal} from \textbf{Theorem} \textit{2}, as,

\begin{align*}
             \mathscr{E} \big[ \mathcal{G}_N \big]  \leq C \big( \mathrm{log} N \big)^{\frac{3}{4}}  \text{ } \text{ . } 
\end{align*}

\noindent For the remaining item, to obtain the desired lower bound on the correlation length from the upper bound provided above on the \textit{greedy lattice animal}, observe, immediately,

\begin{align*}
       \mathscr{L} \geq \mathrm{exp} \big( \mathcal{C}^{-1} \epsilon^{\frac{4}{3}}  \big)     \text{ } \text{ , } 
\end{align*}

\noindent as a result of the fact that, for $\epsilon$ sufficiently small, with high probability,

\begin{align*}
      \textbf{P} \bigg[     \underset{B \in \mathcal{B}_N}{\mathrm{max}}        \frac{\underset{0 \leq \alpha \leq 1-q}{\underset{i \in A}{\mathlarger{\sum}}}   h^{\alpha}_i   }{| \partial B |}      {\sim} \big( \mathrm{log}N \big)^{\frac{3}{4}} \leq \epsilon^{-1}  \bigg] \approx 1 \text{ } \text{ , }
\end{align*}

\noindent for,

\begin{align*}
  N \leq \mathrm{exp} \big( \mathcal{C}^{-1} \epsilon^{\frac{4}{3}} \big) \text{ } \text{ . } 
\end{align*}

\noindent Finally, observe,

\begin{align*}
      \sigma^{\mathrm{w}}_0 \big( \Lambda_N , h \big) \equiv \alpha   \Longleftrightarrow \underset{B \in \mathcal{B}_N}{\mathrm{max}}        \frac{\underset{0 \leq \alpha \leq 1-q}{\underset{i \in A}{\mathlarger{\sum}}}   h^{\alpha}_i   }{| \partial B |}   \geq 1
 \text{ } \text{ , } 
\end{align*}

\noindent for $0 \leq \alpha \leq q-1$, for the spin at the origin, from which we conclude the argument. \boxed{}

\section{References}

\noindent [1] Angles d'Auriac, J.-C., Preissmann, M., Rammal, R. The random field Ising model: algorithmic complexity and phase transition. \textit{Le Journal de Physique-Lettres} \textbf{46} (1985).

\bigskip

\noindent [2] Angles d'Aurias, J.-C., Sourlas, N. The 3-d Random Field Ising Model at zero temperature. \textit{arXiv: 9704088 v1} (1997).

\bigskip

\noindent [3] Cerruti, B., Vives, E. Random Field Potts model with dipolar-like interactions: hysteresis, avalanches, and microstructure. \textit{arXiv: 0710.0286 v1} (2007).

\bigskip

\noindent [4] Ding, J., Song, J., Sun, R. A New Correlation Inequality for Ising Models with External Fields. \textit{arXiv: 2107.09243 v3} (2022).

\bigskip

\noindent [5] Ding, J., Wirth, W. Correlation length of the two-dimensional random field Ising model via greedy lattice animal. \textit{arXiv: 2011.08768 v3} (2022).

\bigskip

\noindent [6] Ding, J., Zhuang, Z. Long range order for the random field Ising and Potts models. \textit{arXiv: 2210.0453 v1} (2021).

\bigskip

\noindent [7] Fytas, N.G., Martin-Mayor, V., Picco, M., Sourlas, N. Review of recent developments in the random-field Ising model. \textit{arXiv: 1711.09597 v3} (2017).

\bigskip

\noindent [8] Fytas, N.G., Martin-Mayor, V. Universality in three-dimensional random-field Ising model. \textit{arXiv: 1304.0318 v2} (2013).

\bigskip

\noindent [9] Goldschmidt, Y.Y., Xu, G. Phase diagram of the random-field Potts model in three dimensions. \textit{Phys. Rev. B} \textbf{32} (1985).

\bigskip

\noindent [10] Grinstein, G., Ma,S-K. Surface tension, roughening, and lower critical dimension in the random-field Ising model. \textit{Phys. Rev. B} \textbf{28}, 2588 (1983).

\bigskip

\noindent [11] Imbrie, J.Z. Lower Critical Dimension of the Random-Field Ising Model. \textit{Phys. Rev. Lett.} \textbf{53}, 1747 (1984).

\bigskip

\noindent [12] Kumar, M., Banerjee, V., Puri, S., Weigel, M. Critical behavior of the three-state random-field Potts model in three dimensions. \textit{arXiv: 2205.13047 v2} (2022).

\bigskip

\noindent [13] Kumar, M., Kumar, R., Weigel, M., Banerjee, V., Janke, W., Puri, S. Approximate ground states of the random-field Potts model from graph ctus. \textit{Phys. Rev. E 97}, 053307 (2018).

\bigskip

\noindent [14] Kumar, M., Weigel, M. Quasi-exact ground-state algorithm for the random-field Potts model. \textit{arXiv:2204.11745 v1} (2022).

\bigskip

\noindent [15] Krzakala, F. Ricci-Tersenghi, F. Zdeborova, L. Elusive Spin Glass Phase in the Random Field Ising Model. \textit{arXiv: 0911.15551 v2} (2010).

\bigskip

\noindent [16] Mezard, M., Young, A.P. Replica Symmetry Breaking in the Random Field Ising Model. \textit{Europhysics Letters}, \textbf{18}, 7 (1992).

\bigskip

\noindent [17] Nattermann, T. Theory of the Random Field Ising Model. \textit{arXiv: 9705295 v1} (1997).

\bigskip

\noindent [18] Qian, H., Bursill, L.A. Random-field Potts model for the Polar Domains of Lead Magnesium Niobate and Lead Scandium Tantalate. \textit{International Journal of of Modern Physics B} \textbf{10} 16 (1996). 

\bigskip

\noindent [19] Reed, P. The Potts model in a random field: a Monte Carlo study. \textit{Journal of Physics c: Solid State Physics}, \textbf{18}, 20 (1985).

\bigskip

\noindent [20] Son, S-W., Jeong, H. Noh, J.D. Random field Ising model and community structure in complex networks. \textit{arXiv: 0502672 v1} (2005).

\bigskip

\noindent [21] Talagrand, M. Regularity of gaussian processes. \textit{acta Mathematica} \textbf{159}, 99-149 (1987).

\bigskip

\noindent [22] Talagrand, M. Upper and Lower Bounds for Stochastic Processes. \textit{Springer-Verlag} (2021).

\bigskip
 
\noindent [23] Tissier, M., Tarjus, G. Supersymmetry and its spontaneous breaking in the random field Ising model. \textit{arXiv: 1103.4812 v2} (2011).

\bigskip

\noindent [24] Turkoglu, A., Berker, A.N. Phase transitions of the variety of the random-fields Potts models. \textit{Physica A: Statistical Mechanics and its Applications} \textbf{583} 126339 (2021).

\bigskip

\noindent [25] Wu, Y., Machta, J. Ground states and thermal states of the random field Ising model. \textit{arXiv: 0501619 v2} (2005).

\end{document}